\newcommand{\ncm}{\newcommand}
\ncm{\cstar}{$C^{*}$-algebra} \ncm{\cstars}{$C^{*}$-algebras} 
\ncm{\Ga}{\mbox{$\Gamma$}} \ncm{\ra}{\rightarrow} 
\ncm{\st}[2]{\stackrel{#1}{#2}} \ncm{\ral}{\longrightarrow} 
\ncm{\al}{\mbox{$\alpha $}} \ncm{\gam}{\mbox{$\gamma $}} 
\ncm{\vp}{\varphi} \ncm{\vep}{\varepsilon} \ncm{\vt}{\vartheta} 
\ncm{\Ad}{\mbox{\rm Ad}} \ncm{\mod}{\mbox{\rm mod}} 
\ncm{\ol}{\overline} \ncm{\OA}[1]{{\cal O}_{#1}} \ncm{\OPA}[1] 
{\mbox{$\bar{\cal O}_{#1}$}} \ncm{\FA}[1]{{\cal F}_{#1}} 
\ncm{\SFA}[1] {\mbox{$\bar{\cal F}_{#1}$}} \ncm{\EA}[1]{{\cal 
E}_{#1}} \ncm{\KO}[2]{K_{#1} (\OA{#2})} \ncm{\KE}[2]{K_{#1} 
(\EA{#2})} \ncm{\ld}{\lneqq} \ncm{\lb}{\label} \ncm{\qed}{\hfill 
$\square$ \smallskip} 

\documentclass{article}[10pt]
\usepackage{latexsym,amssymb}

\newtheorem{theo}{Theorem}[section]
\newtheorem{cor}[theo]{Corollary}
\newtheorem{notation}[theo]{Notation}
\newtheorem{lem}[theo]{Lemma}
\newtheorem{prop}[theo]{Proposition}
\newtheorem{remark}[theo]{Remark}
\newtheorem{example}[theo]{Example}
\newtheorem{definition}[theo]{Definition}
\newtheorem{problem}[theo]{Problem}

\newenvironment{rem}{\begin{remark}\rm}{\end{remark}}
\newenvironment{ex}{\begin{example}\rm}{\end{example}}

\newenvironment{pf}{{\it Proof:}\pagebreak[1]}
{\qed} 
\newenvironment{df}{\begin{definition}\rm}{\end{definition}}

\newenvironment{jbas}{\begin{equation}}{\end{equation}}
\newenvironment{jbam}{\begin{eqnarray*}}{\end{eqnarray*}}

\newenvironment{jba}{\begin{displaymath}}{\end{displaymath}}
\ncm{\R}{\mbox{\bf R}} 

\ncm{\Z}{\mathbb Z} 

\ncm{\T}{\mathbb T} 

\ncm{\TT}{\T$^{2}$} 

\ncm{\N}{\mathbb N} 

\ncm{\C}{\mathbb C} 

\begin{document}
\begin{center} {\Large  MAXIMAL ABELIAN SUBALGEBRAS OF ${\cal O}_n$ \\}
\end{center} \bigskip \begin{center} {E. J. Beggs and P. Goldstein \\
24 October 2000}
\end{center}
\bigskip
\subsection*{Abstract}
We consider maximal abelian subalgebras of $\OA{n}$ which are 
invariant to the standard circle action. It turns out that these 
are all contained in the zero grade of $\OA{n}$. Then we consider 
shift invariant maximal abelian subalgebras of the zero grade, 
which are also invariant to a ``second shift'' map, and show that 
these are just infinite tensor products of diagonal matrices in 
the standard UHF picture of the zero grade.

\section{Introduction}

In this note, we are concerned with certain abelian subalgebras of 
Cuntz algebra $\OA{n}$. Let $\OA{n} = C^{*} (s_1 , \ldots, s_n )$. 
As usual, for $\mu = i_{1}\ldots i_{k}$, $i_{j} \in \{ 1, \ldots, 
n \}$, we let $|\mu|=k$ be the length of $\mu$ and denote 
$s_{i_{1}}\ldots s_{i_{k}}$ by $s_\mu$. The set of all finite 
words in $\{ 1, \ldots, n \}$ is denoted by ${\cal W}(n)$. Let 
\begin{jbas}
\sigma (x) = \sum_{i=1}^{n} s_i x s_{i}^{*} 
\end{jbas}
be the canonical endomorphism on $\OA{n}$, 
\begin{jbas}
\omega_t (s_i) = t s_i , \ i=1,\ldots, n, \ t \in \T
\end{jbas}the standard circle action, and consider the $C^*$--subalgebra
$\cal D$ of $\OA{n}$ defined as 
\begin{jbas}
{\cal D} = C^* \{ s_{\mu} s_{\mu}^* ; \mu \in {\cal W}(n) \} 
\end{jbas}Then $\cal D$ is an abelian subalgebra (cf. \cite{cu1}),
and below we list some of its properties: 
\begin{itemize}
\item $\sigma({\cal D}) \subset {\cal D}$
\item $\omega_t ({\cal D}) \subset {\cal D}$; in fact, even more is
true: $\omega_t (d) =d $, for all $t \in \T$ and $d \in {\cal D}$; 
hence, ${\cal D} \subset \OA{n}^{0}$, the fixed--point algebra 
under the circle action 
\item $\cal D$ is a maximal abelian subalgebra of $\OA{n}^0$;
furthermore, $\cal D$ is maximal abelian in $\OA{n}$ (cf. 
\cite{cuk}) 
\end{itemize}

The subject of the present work is to give a characterisation of 
$\cal D$ in the above terms. More precisely, we prove (cf. 
\ref{mth1}): 

\begin{theo}
Let $A$ be a maximal abelian subalgebra of $\OA{n}$ such that 
$\omega (A) \subset A$, $\sigma (A) \subset A$ and $\tilde{\sigma} 
(A) \subset A$. Then there is an automorphism $\alpha_U$ of 
$\OA{n}$, determined by a unitary $U \in M_n(\C)$, such that 
$\alpha_U (A) = {\cal D}$. Furthermore, $\alpha_U$ commutes with 
$\omega$, $\sigma$ and $\tilde{\sigma}$. 
\end{theo}

We now describe the above notation. Let $\psi : \OA{n} \ra M_n 
\otimes \OA{n}$ be given by 
\begin{jbas}
x \stackrel{\psi}{\mapsto} \left[ \begin{array}{ccc} s_{1}^{*} x 
s_{1} & \ldots & s_{1}^{*} x s_{n} \\ \vdots & & \vdots \\ 
s_{n}^{*} x s_{1} & \ldots & s_{n}^{*} x s_{n}  \end{array} 
\right] 
\end{jbas}Then $\psi$ is an isomorphism from $\OA{n}$ onto $M_n \otimes
\OA{n}$ (cf. \cite{cho}). 

The endomorphism $\tilde{\sigma} : \OA{n} \ra \OA{n}$  is the 
``second shift'', given by the formula 
\begin{jbas}
\tilde{\sigma} = \sum_{i,j,k} s_i s_k s_{i}^{*} x s_j s_{k}^{*} 
s_{j}^{*}, 
\end{jbas}and determined by the
following diagram 
\begin{jbas}
\begin{array}{ccc}
\OA{n} & \stackrel{\psi}{\ral} & M_n \otimes \OA{n} \\ 
\tilde{\sigma} \downarrow & & id \downarrow \otimes \sigma  \\ 
\OA{n} & \stackrel{\psi^{-1}}{\longleftarrow} & M_n \otimes \OA{n} 
\end{array}
\end{jbas}Note that it follows easily that
\begin{jba}
\tilde{\sigma}(1) = 1, \ \tilde{\sigma}(s_i s_{j}^{*}) = s_i 
s_{j}^{*}, 
\end{jba}and
\begin{jba}
\tilde{\sigma}( s_i s_{\mu} s_{\nu}^{*} s_{j}^{*}) = s_i \sigma( 
s_{\mu} s_{\nu}^{*} ) s_{j}^{*}, 
\end{jba}for all $\mu , \nu \in {\cal W}(n)$ such that $|\mu| =
|\nu|$. 

Finally, for $u = [u_{ij}]_{i,j=1}^{n}$ a unitary in $M_{n}(\C)$, 
let $U= \sum_{i,j} u_{ij}s_{i}s_{j}^{*}$. Then $U$ is a unitary in 
$\OA{n}$, and the map $s_i \mapsto U s_i , \ i=1,\ldots,n$ extends 
to an isomorphism of $\OA{n}$, denoted $\alpha_U$. 

We also prove the same result as \ref{mth1}, but starting from 
slightly different assumptions (cf. \ref{mth2}): 

\begin{theo}
Let $A$ be a maximal abelian algebra in $\OA{n}$ such that $A \cap 
\OA{n}^{0}$ is maximal abelian in $\OA{n}^{0}$, $\sigma (A) 
\subset A$ and $\tilde{\sigma} (A) \subset A$. Then there is an 
automorphism $\alpha_U$ of $\OA{n}$, determined by a unitary $U 
\in M_n(\C)$, such that $\alpha_U (A) = {\cal D}$. Furthermore, 
$\alpha_U$ commutes with $\sigma$ and $\tilde{\sigma}$. 
\end{theo}

The paper is organised as follows. In Section 2, we show that an 
algebra which is maximal in the class of abelian algebras that are 
invariant under the action of the circle is indeed maximal 
abelian. In Section 3, we describe maximal abelian subalgebras 
that are invariant under two shift maps. Both Section 2 and 3 are 
done in a slightly more general setting. Finally, in Section 4, we 
apply these results to the particular case of $\OA{n}$ and easily 
obtain the stated characterisation of the abelian subalgebra $A$.

\subsection*{Acknowledgements} The authors would like to thank G. 
Elliott, R. Exel, S. Stratila, S. Wassermann and N-C. Wong for 
very useful comments. This research was supported by an EPSRC 
Research Assistanship (P.G.). 

\section{Maximal (abelian $\T$--invariant *--subalgebras)}

Take a $C^*$--algebra $B$, and let $\omega : \T \ral Aut(B)$ be a 
homomorphism that is continuous in the topology of pointwise 
convergence. This means that for each $b \in B$ the map $t \mapsto 
\omega_{t} (b)$ is continuous, and the triple $(B, \T, \omega)$ is 
called a $C^*$--dynamical system (cf. \cite[7.4.1]{ped}).  

Consider the class of *--subalgebras $A$ of $B$ which are abelian 
and $\T$--invariant (i.e. $\omega_{t} (a)\in A$ for all $a\in A$ 
and all $t \in \T$). We will call a subalgebra which is maximal in 
this class a maximal (abelian $\T$-invariant *--subalgebra), using 
brackets to avoid ambiguity. Our task is to show that we can 
remove the brackets, that is, we show that a maximal (abelian 
$\T$-invariant *--subalgebra) is actually maximal abelian. 

The main result of this section is \ref{th1}. For convenience we 
assume that $B$ is a subalgebra of $B(H)$ for some Hilbert space 
$H$. 

\begin{df}
Let $(B, \T, \omega)$ be a $C^*$--dynamical system, and define 
$B_n = \{ b \in B: \ \omega_t (b) = t^n b , \mbox { for all } t 
\in \T \}$. Let $\pi_n : B \ra B$ be defined as 
\begin{jba}
\pi_n (b) = \int_{\T} t^{-n} \omega_t (b) dt,
\end{jba}where $dt$ is the Haar measure on $\T$ (i.e. normalised 
Lebesgue measure). Then each $B_n$ is a closed linear subspace in 
$B$, and each $\pi_n$ is a linear contraction with image $B_n$. 
Also $B_0$ is a subalgebra and $\pi_0$ is a conditional 
expectation to $B_0$. Furthermore, we have 
\begin{jba} 
\pi_m (b) = \delta_{n,m} b, \ b \in B_n.
\end{jba}
\end{df}
 
\begin{prop} The commutant $A'\cap B$ of $A$ in $B$ is an
$\T$--invariant *-subalgebra of $B$, and $A\subset A'$. Further 
the image of $A'\cap B$ under $\pi_n$ is contained in $A'\cap B$. 
\end{prop}
\begin{pf}
Suppose that $b\in A'\cap B$. Since $A$ is a *-subalgebra, for all 
$a\in A$, $a^*b=ba^*$, so we see that $b^*a=ab^*$. Likewise, for 
any $t \in \T$, $\omega_{t^{-1}}(a)b=b \omega_{t^{-1}}(a)$ so 
$a\omega_t (b)=\omega_t (b) a$. As $A$ is abelian, $A\subset A'$. 
Finally we take the equation $a\omega_t (b)=\omega_t (b) a$, 
divide by $t^n$ and integrate to see that $\pi_n (b)\in A'$.  
\end{pf}

\begin{prop} The image of $A'\cap B$
under $\pi_0$ is contained in $A$. 
\end{prop}
\begin{pf}
Suppose that $b\in A'\cap B$. By considering $b+b^*$ and 
$i(b-b^*)$ we may suppose that $b$ is actually Hermitian. Then 
$\pi_0(b)$ is Hermitian and fixed by the circle action, so the 
algebra generated by $A$ and $\pi_0(b)$ is an abelian circle 
invariant *-subalgebra of $B$, and so $\pi_0(b)\in A$ by 
maximality.
\end{pf}

\begin{prop} The image of $A'\cap B$
under $\pi_n$ is contained in $A$ for all $n\in\mathbb Z$. 
\end{prop}
\begin{pf}
Suppose that $b\in A'\cap B$. Then $\pi_n(b)\in A'$, and 
$\pi_n(b)\pi_n(b)^*$ is an Hermitian circle invariant element of 
$A'\cap B$. By maximality we then have $\pi_n(b)\pi_n(b)^*\in A$, 
and similarly we have $\pi_n(b)^*\pi_n(b)\in A$. This means that 
$\pi_n(b)$ commutes with both $\pi_n(b)\pi_n(b)^*$ and 
$\pi_n(b)^*\pi_n(b)$, and so is normal by the next lemma. Now the 
algebra generated by $A$, $\pi_n(b)$ and $\pi_n(b)^*$ is an 
abelian circle invariant *-algebra, so $\pi_n(b)\in A$ by 
maximality.
\end{pf}

We initially proved the next lemma using polar decomposition. The 
following much simpler proof is due independently to S. Wassermann 
and N--C. Wong: 

\begin{lem}
Let $x \in B$ commute with both $x x^*$ and $x^* x$. Then $x$ is 
normal. 
\end{lem}
\begin{pf}
We need to show that $x x^* - x^* x =0$. Since $x x^* - x^* x$ is 
selfadjoint, that is equivalent to $( x x^* - x^* x )^2 = 0$, 
which follows immediately, since the assumption implies 
\begin{jba}
x^* x x x^* = x (x^* x) x^* \mbox{ and } x x^* x^* x = x^* (x x^* 
) x 
\end{jba}
\end{pf}

\begin{prop} $A'\cap B\subset A''$.
\end{prop}
\begin{pf}
Take $b\in A'\cap B$. For any $\xi,\eta\in H$ we define a 
continuous function $f:\T \ra \mathbb C$ by 
$f(t)=\big<\xi,\omega_t (b)(\eta)\big>$. By Fourier analysis we 
get Fourier coefficients $f_n=\big<\xi,\pi_n(b)(\eta)\big>$, where 
\begin{jbas}
\sum_{n=-m}^m t^n\, f_n\ \ral f(t), \  t \in \T 
\end{jbas}in the  $L^2(\T)$ topology as $m\to\infty$. Since $B \subset 
B(H)$, we can put $\eta=c(\kappa)$ for some $c\in A'$ and 
$\kappa\in H$. Then as $\pi_n(b)\in A$ we see that 
$f_n=\big<\xi,c\pi_n(b)(\kappa)\big>=\big<c^*\xi,\pi_n(b)(\kappa)\big>$. 
Now we can write 
\begin{jba}
\sum_{n=-m}^m t^n  f_n  \ra  \big<c^*\xi,\omega(b)(\kappa)\big>\ 
=\ \big<\xi,c\omega(b)(\kappa)\big>, \  t \in \T
\end{jba}in the $L^2(\T)$ topology as $m\to\infty$. The two limits are the 
same in $L^2(\T)$, so $\big<\xi,c\omega_t 
(b)(\kappa)\big>=\big<\xi,\omega_t (b)c(\kappa)\big>$ almost 
everywhere in $\T$. By continuity they are the same at $t=1$, so 
$cb=bc$.
\end{pf}

\begin{theo}
\lb{th1} Let $(B, \T, \omega)$ be a $C^*$--dynamical system, and 
suppose that $A$ is a maximal (abelian $\T$--invariant 
*--subalgebra) of $B$. Then $A$ is a maximal abelian subalgebra of 
$B$.
\end{theo} 
\begin{pf}
By the previous proposition we see that $A'\cap B$ is abelian. Now 
$A'\cap B$ is an abelian circle invariant *-subalgebra of $B$ 
which contains $A$, so $A=A'\cap B$ by the maximality condition on 
$A$. But the equation $A=A'\cap B$ means that $A$ is maximal 
abelian in $B$. 
\end{pf}

The following example -- due to R. Exel -- shows that the previous 
theorem does not hold for an arbitrary dynamical system  $(B, G, 
\omega)$, even with $G$ compact:  
\begin{ex} 
Consider the adjoint action of $SU_2$ on $M_2 (\C)$. The 
subalgebra consisting of the complex multiples of the identity is 
maximal among the class of abelian $SU_2$--invariant 
$*$--subalgebras. However, it is not maximal abelian, as it is 
properly contained in the diagonal matrices.  
\end{ex}

\section{Maximal abelian *-subalgebras of ${\cal O}_n$
contained in the zero grade} 

Let $B$ be a unital $C^*$--algebra with a given isomorphism 
$\psi:B\to M_n\otimes B$ with $\psi(1)=I_n\otimes 1$, $I_n$ being 
the identity matrix in $M_n$. We define isomorphisms $\psi_m:B\to 
(M_n)^{\otimes m}\otimes B$ ($m\ge 0$) recursively, beginning with 
$\psi_0:B\to B$ the identity, $\psi_1=\psi$, and continuing by 
defining $\psi_{m+1}$ to be the composition \[ B  
\stackrel{\psi_m} \longrightarrow (M_n)^{\otimes m}\otimes B 
\stackrel{{\rm id}^{\otimes m}\otimes \psi} 
  \longrightarrow
(M_n)^{\otimes m+1}\otimes B\ , 
\]
where ${\rm id}:M_n\to M_n$ is the identity map. Now we define an 
algebra map $\kappa_m:M_n^{\otimes m}\to B$ by 
$\kappa_m(x)=\psi_m^{-1}(x\otimes 1)$. Since $\psi(1)=I_n\otimes 
1$ we get the commutative diagram \begin{jba} 
\begin{array}{ccc} M_n^{\otimes m} & \stackrel{\kappa_m}{\ral} & B \\
\quad \downarrow {\rm id}\otimes I_n & & \quad \downarrow {\rm 
id}_B \\ M_n^{\otimes m+1} & \stackrel{\kappa_{m+1}}{\ral} & B 
\end{array} 
\end{jba} Define $C\subset B$ to be the closure of the union of the
subalgebras $\kappa_m(M_n^{\otimes m})$. 

We can define shift maps $\sigma_m:B\to B$ ($m\ge 1$) by the 
composition 
\[
B  \stackrel{\psi_{m-1}} \longrightarrow (M_n)^{\otimes 
m-1}\otimes B \stackrel{{\rm id}^{\otimes m-1}\otimes f} 
  \longrightarrow
(M_n)^{\otimes m}\otimes B \stackrel{\psi_{m}^{-1}} 
  \longrightarrow B\ ,
\]
where $f:B\to M_n\otimes B$ is the algebra map $f(b)=I_n\otimes 
b$. 

Let $E_{ij}\in M_n$ be the matrix with entry $1$ in row $i$ column 
$j$, and zeros elsewhere. Define a linear map 
$e_{ij}:M_n\to\mathbb C$ by $e_{ij}(E_{kl}) 
=\delta_{ik}\delta_{jl}$. Now we can define a map $\chi_{mij}:B\to 
B$ ($m\ge 1$) by the composition \[ B \stackrel{\psi_m} 
\longrightarrow (M_n)^{\otimes m}\otimes B 
  \stackrel{{\rm id}^{\otimes m-1}\otimes
e_{ij}\otimes {\rm id_B}} \longrightarrow (M_n)^{\otimes 
m-1}\otimes B 
  \stackrel{\psi_{m-1}^{-1}} \longrightarrow B\ .
\]

\begin{prop} For all $b\in B$ and $y\in M_n\otimes B$,
$(e_{ij}\otimes {\rm id}_B)(f(b).y)=b.((e_{ij}\otimes {\rm 
id}_B)(y))$ and $(e_{ij}\otimes {\rm 
id}_B)(y.f(b))=((e_{ij}\otimes {\rm id}_B)(y)).b$. 
\end{prop} 
\begin{pf}
Take $y=y_1\otimes y_2\in M_n\otimes B$ (linear combinations 
of terms of this form are dense in $M_n\otimes B$). Then 
\begin{jbam}
(e_{ij}\otimes {\rm id}_B)(f(b).y)\ =\ (e_{ij}\otimes {\rm 
id}_B)((I\otimes b)(y_1\otimes y_2))\ = \\ (e_{ij}\otimes {\rm 
id}_B)(y_1\otimes by_2)\ =\ e_{ij}(y_1)\, by_2 \ =\ 
b.((e_{ij}\otimes {\rm id}_B)(y))\ 
\end{jbam}The other way round is the same. 
\end{pf}

\begin{cor} For all $b,c\in B$,
$\chi_{mij}(\sigma_m(b).c)=b.\chi_{mij}(c)$ and 
$\chi_{mij}(c.\sigma_m(b))=\chi_{mij}(c).b$. 
\end{cor}
\begin{pf}
This is essentially the same as the previous proposition.
\end{pf}

\begin{cor} Suppose that $A$ is a maximal abelian $*$-subalgebra of $B$,
obeying the condition $\sigma_m(A)\subset A$. Then for all $1\le 
i,j\le n$, $\chi_{mij}(A)\subset A$. 
\end{cor}
\begin{pf}
Take $a\in A$. Then for all $a'\in A$ we have 
$\sigma_m(a').a=a.\sigma_m(a')$. Applying $\chi_{mij}$ to this we 
get $a'.\chi_{mij}(a)=\chi_{mij}(a).a'$, so $\chi_{mij}(a)\in A$ 
by maximality.
\end{pf}

\begin{prop} Suppose that $A$ is a maximal abelian $*$-subalgebra of $B$,
obeying the condition $\sigma_1(A)\subset A$. Then 
$\psi_m(A)\subset M_n^{\otimes m}\otimes A$. 
\end{prop} 
\begin{pf}
First note that $\psi(a)=\sum_{ij}E_{ij}\otimes\chi_{1ij}(a)$, so 
$\psi(A)\subset M_n\otimes A$ by the last proposition. The rest 
follows by induction. 
\end{pf}

\begin{df}
Take a unital algebra map $\phi:A\to \mathbb C$, and extend it to 
a positive contraction $\phi:B\to \mathbb C$. Then we define a map 
$\phi_m:B\to M_n^{\otimes m}$ by 
\[
B  \stackrel{\psi_m} \longrightarrow M_n^{\otimes m}\otimes B 
  \stackrel{{\rm id}^{\otimes m}\otimes \phi} \longrightarrow
M_n^{\otimes m}\ . 
\]
Since $\phi_m(1)=1$, it follows from \cite[3.1.6]{ped} that 
$\phi_m$ is a contraction. On the other hand, this is clearly a 
unital homomorphism when restricted to $A$. We denote by $D$ the 
image of $\phi_1:A\to M_n$. 
\end{df}

\begin{prop} If $\sigma_1 (A)\subset A$, then
$\phi_{m+1}(A)\subset M_n\otimes \phi_{m}(A)$. 
\end{prop}
\begin{pf}
 We can write $\phi_{m+1}$ as 
\[
A  \stackrel{\psi_1} \longrightarrow M_n\otimes A  \stackrel{{\rm 
id}\otimes\psi_m} \longrightarrow (M_n)^{\otimes m+1}\otimes A 
  \stackrel{{\rm id}^{\otimes m+1}\otimes \phi} \longrightarrow
M_n^{\otimes m+1}\ , 
\]
which can be rewritten as 
\[
A  \stackrel{\psi_1} \longrightarrow M_n\otimes A  \stackrel{{\rm 
id}\otimes\phi_m} \longrightarrow M_n^{\otimes m+1}\ , 
\]
so we see that $\phi_{m+1}(A)\subset M_n\otimes \phi_{m}(A)$. 
\end{pf}

\begin{prop} If $\sigma_2 (A)\subset A$, then
$({\rm id}\otimes e_{ij}\otimes {\rm id}^{\otimes m-1}) 
\phi_{m+1}(A)\subset \phi_{m}(A)$ for $m\ge 1$. 
\end{prop}
\begin{pf}
The map $({\rm id}\otimes e_{ij} \otimes {\rm id}^{\otimes 
m-1})\circ \phi_{m+1}$ is 
\[
B  \stackrel{\psi_{m+1}} \longrightarrow M_n^{\otimes m+1}\otimes 
B 
   \stackrel{{\rm id}^{\otimes m+1}\otimes\phi} \longrightarrow
M_n^{\otimes m+1} 
  \stackrel{{\rm id}\otimes e_{ij}
\otimes {\rm id}^{\otimes m-1}} \longrightarrow M_n^{\otimes m}\ , 
\]
which can be rewritten as 
\[
B  \stackrel{\psi_{m+1}} \longrightarrow M_n^{\otimes m+1}\otimes 
B 
  \stackrel{{\rm id}\otimes e_{ij}
\otimes {\rm id}^{\otimes m-1}\otimes{\rm id}_B} \longrightarrow 
M_n^{\otimes m}\otimes B \stackrel{\psi_m^{-1}} \longrightarrow B 
\stackrel{\psi_m} \longrightarrow M_n^{\otimes m}\otimes B 
   \stackrel{{\rm id}^{\otimes m}\otimes\phi} \longrightarrow
M_n^{\otimes m}\ . 
\]
This can be shown to be 
\[
B  \stackrel{\psi_{2}} \longrightarrow M_n^{\otimes 2}\otimes B 
  \stackrel{{\rm id}\otimes e_{ij}
\otimes{\rm id}_B} \longrightarrow M_n\otimes B 
\stackrel{\psi_1^{-1}} \longrightarrow B \stackrel{\psi_m} 
\longrightarrow M_n^{\otimes m}\otimes B 
   \stackrel{{\rm id}^{\otimes m}\otimes\phi} \longrightarrow
M_n^{\otimes m}\ , 
\]
so we see that $({\rm id}\otimes e_{ij} \otimes {\rm id}^{\otimes 
m-1})\circ \phi_{m+1}=\phi_m\circ \chi_{2ij}:B\to M_n^{\otimes 
m}$. Now use $\chi_{2ij}(A)\subset A$. 
\end{pf}

\begin{cor}
\lb{c3.8}
If $\sigma_1 (A)\subset A$ and $\sigma_2 (A)\subset A$, 
then $\phi_m(A)\subset D^{\otimes m}$. 
\end{cor}
\begin{pf} 
This is proved by induction. First note that $\phi_1(A)\subset 
D^{\otimes 1}$ by definition of $D$. Now assume that 
$\phi_m(A)\subset D^{\otimes m}$, and consider $m+1$. By the 
previous proposition we see that $\phi_{m+1}(A)\subset D\otimes 
M_n\otimes D^{\otimes m-1}$, whereas the proposition before that 
says that $\phi_{m+1}(A)\subset M_n\otimes D^{\otimes m}$. The 
result follows by standard linear algebra. 
\end{pf}

\begin{prop}
\lb{p3.9}
Given $c\in C$ and $\epsilon>0$, there is an $m\ge 1$ so 
that $|\kappa_m(\phi_m(c))-c|<\epsilon$. 
\end{prop}
\begin{pf}
There is an $m\ge 1$ and an $x\in M_n^{\otimes m}$ so that 
$|c-\kappa_m(x)|<\epsilon/2$. Since $\phi(1)=1$ we get 
$\phi_m(\kappa_m(x))=x$, and since $\phi_m$ is a contraction, 
$|\phi_m(c)-x|<\epsilon/2$. Finally as $\kappa_m$ is a 
contraction, $|\kappa_m(\phi_m(c))-\kappa_m(x)|<\epsilon/2$. 
\end{pf} 

Now, let $D^{\infty}$ stand for the closure of the union of 
$\kappa_m (D^{\otimes m})$ for $m \geq 1$. Then we have: 

\begin{theo}
\lb{mth4} Suppose that $A \cap C$ is maximal abelian in $C$. Then 
$A \cap C = D^{\infty}$ and $D$ is maximal abelian in $M_n (\C)$. 
\end{theo}
\begin{pf}
\ref{c3.8} and \ref{p3.9} show that $A \cap C \subset D^{\infty}$. 
Since $D^{\infty}$ is abelian and $A \cap C$ is maximal, it 
follows that $A \cap C=D^{\infty}$. Hence, $D$ is maximal abelian 
in $M_n(\C)$. 
\end{pf}

\section{Maximal (abelian $\T$--invariant *-subalgebras) of ${\cal O}_n$}

In this section, we apply the results from Section 2 and 3 to 
maximal abelian subalgebras of $\OA{n}$ that are invariant under 
the standard circle action. The notation is as in the 
Introduction. 

The next lemma is probably well--known, but we couldn't find a 
reference: 
\begin{lem}
Let $x$ be in $\OA{n}^k$ (i.e. $\omega_t (x)= t^k x$), for $k \neq 
0$. If $x$ is normal , then $x = 0$. 
\end{lem}
\begin{pf}
Suppose $k > 0$. Let $y= x (s_{1}^*)^k \in \OA{n}^0$, and let 
$\tau$ be the faithful normalised trace on $\OA{n}^0 \cong 
M_{n^\infty}(\C)$. Then $y y^* = x x^*$, $y^* y = s_{1}^k x^* x 
(s_{1}^* )^k$, and $\tau ( y y^* ) = \tau (y^* y )$ imply 
\begin{jba}
\tau (x x^* ) = n^{-k} \tau ( x^* x)
\end{jba}If $x x^* = x^* x$, then
\begin{jba}
\tau ( x^* x) = n^{-k} \tau ( x^* x), 
\end{jba}hence $\tau (x^* x) = 0$.
\end{pf}
\begin{theo}
\lb{mth3} Suppose that  $A$ is a maximal ( abelian $\T$--invariant 
*-subalgebra) of ${\cal O}_n$. Then $A \subset \OA{n}^0$. 
\end{theo}
\begin{pf}
By the previous lemma, $\pi_k(a)=0$ for all $a\in A$ and $k\neq 
0$. By Fourier analysis we get 
\begin{jba}
\sum_{k=-m}^m t^k  \big<\xi,\pi_k(a)(\eta)\big> \ \to\ 
\big<\xi,\omega_t (a)(\eta)\big>, \ t \in \T 
\end{jba}in the  $L^2(\T)$ topology as $m\to\infty$. But then 
$\big<\xi,\omega_t (a)(\eta)\big>$ is constant on $\T$, so 
$\omega_t (a)=a$.
\end{pf}
\begin{rem}
Note that, if $B$ is a $C^*$--algebra of the form $B = A 
\rtimes_\alpha \N$, where $A$ has a normalised faithful trace,  
$\alpha$ is a trace--scaling endomorphism, and the circle action 
is just the dual action with respect to this crossed--product 
representation, the above argument shows that any maximal (abelian 
$\T$--invariant *-subalgebra) will be contained in the 
fixed--point algebra for the circle action. This will be the case 
if $B$ is a simple Cuntz--Krieger algebra, or more generally, for 
certain Cuntz--Pimsner algebras (cf. \cite{pim}). 
\end{rem}
\begin{theo}
\lb{mth1} 
Let $A$ be a maximal abelian subalgebra of $\OA{n}$ such 
that $\omega (A) \subset A$, $\sigma (A) \subset A$ and 
$\tilde{\sigma} (A) \subset A$. Then there is an automorphism 
$\alpha_U$ of $\OA{n}$, determined by a unitary $U \in M_n(\C)$, 
such that $\alpha_U (A) = {\cal D}$. Furthermore, $\alpha_U$ 
commutes with $\omega$, $\sigma$ and $\tilde{\sigma}$. 
\end{theo}
\begin{pf}
By \ref{mth3}, $A \subset \OA{n}^0$. The result then follows from 
\ref{mth4}, with $B= \OA{n}^0$, $\sigma = \sigma_1$ and  
$\tilde{\sigma} = \sigma_2$, while $u$ is any unitary in $M_n(\C)$ 
that diagonalises the maximal abelian subalgebra $D$. 
\end{pf}

\begin{theo}
\lb{mth2} Let $A$ be a maximal abelian algebra in $\OA{n}$ such 
that $A \cap \OA{n}^{0}$ is maximal abelian in $\OA{n}^{0}$, 
$\sigma (A) \subset A$ and $\tilde{\sigma} (A) \subset A$. Then 
there is an automorphism $\alpha_U$ of $\OA{n}$, determined by a 
unitary $U \in M_n(\C)$, such that $\alpha_U (A) = {\cal D}$. 
Furthermore, $\alpha_U$ commutes with $\sigma$ and 
$\tilde{\sigma}$. 
\end{theo}
\begin{pf}
We apply \ref{mth4}, with $B=\OA{n}$, and $\sigma$ and 
$\tilde{\sigma}$ as in the previous theorem. That shows that 
$D^\infty \subset A$. Since $D^\infty$ is maximal abelian in 
$\OA{n}$ (cf. \cite[2.18]{cuk}), $A = D^\infty$. Note that this 
means that all of $A$ is contained in the zero grade, although no 
assumptions on the circle action were made. 
\end{pf}

\noindent E. J. Beggs,  Department of Mathematics, University of 
Wales Swansea, \\ Swansea SA2 8PP \\ and \\ P. Goldstein, School 
of Mathematics, Cardiff University, 
\\ 
 Cardiff CF4 4YH 
 
\end{document}